\documentstyle[11pt]{article}

\begin{document}

\author{G. Boyadzhiev}
\title{Comparison principle for non - cooperative elliptic systems}
\date{28.02.2007}
\maketitle

\section{Introduction}
In this paper are considered weakly coupled linear elliptic systems
of the form

\smallskip\ \\(1)\qquad $%
L_Mu=0$ in a bounded domain $\Omega \in R^n$ with smooth boundary

\smallskip\ \\and boundary data $u(x)=g(x)$ on $\partial \Omega $,
where $L_M=L+M$, $L$ is a matrix operator with null off-diagonal
elements $L=diag\left(L_1,  L_2, ... L_N\right) $, and matrix
$M=\{m_{ik}(x)\}_{i,k=1}^{N}$. Scalar operators

\smallskip\ \\ \qquad $L_ku_k = -\sum_{i,j=1}^{n}D_j
\left( a_k^{ij}(x)D_iu_k \right) +\sum_{i=1}^{n}b_k^i(x)D_iu_k+c_ku_k$ in $%
\Omega $

\smallskip\ \\are uniformly elliptic ones for $k=1,2,...N$, i.e. there are constants
$\lambda,\Lambda >0$ such that

\smallskip\ \\(2) \qquad $\lambda \left| \xi \right|
^2\leq \sum_{i,j=1}^{n}a_k^{ij}(x) \xi _i\xi _j\leq \Lambda \left|
\xi \right| ^2$

\smallskip\ \\for every $k$ and any $\xi =(\xi _1,...\xi _n)\in R^n$.

Coefficients $c_k$ and $m_{ik}$ in (1) are supposed continuous in
$\overline{\Omega}$, and $a_k^{ij}(x),b_k^i(x)\in
W^{1,\infty}(\Omega)\cap C(\overline{\Omega})$.

Quasi-linear weakly coupled elliptic systems

\smallskip\ \\(3) \qquad $Q^l(u)=-diva^l(x,u^l,Du^l)+F^l(x,u^1,...u^N,Du^l)=f^l(x)$
in $\Omega $

\smallskip\ \\(4) \qquad $u^l(x)=g^l(x)$ on $\partial \Omega $

\smallskip\ \\$l=1,...N$ are considered as well.

System (3) is supposed  uniformly elliptic one, i.e. there are
continuous and positive functions $\lambda (\left| u\right|
),\Lambda (\left| u\right| )$, $\left| u\right| =\left( \left(
u^1\right) ^2+...+\left( u^N\right) ^2\right) ^{1/2}$, such that
$\lambda (s)$ is monotone-decreasing one, $\Lambda (s)$ is monotone
increasing one and

\smallskip\ \\ \qquad $\lambda (\left| u\right| )\left| \xi ^l\right|
^2\leq \sum_{i,j=1}^{n}\frac{\partial a^{li}}{%
\partial p_j^l}(x,t,u^1,...u^N,p^l)\xi _i^l\xi _j^l\leq \Lambda (\left|
u\right| )\left| \xi ^l\right| ^2$

\smallskip\ \\for every $u^l$ and $\xi ^l=(\xi _1^l,...\xi _n^l)\in R^n$, $%
l=1,2,...N.$

The coefficients $a^l(x,u,p)$, $%
F^l(x,u,p)$, $f^l(x)$, $g^l(x)$ are supposed to be at least
measurable functions with respect to the $x$ variable and locally
Lipschitz continuous on $u^l,u$ and $p$, i.e.

\smallskip\ \\ $
\begin{array}{l}
\left| F^l(x,u,p)-F^l(x,v,q)\right| \leq C(K)\left( \left|
u-v\right|
+\left| p-q\right| \right) , \\
\\
\left| a^l(x,u^l,p)-a^l(x,v^l,q)\right| \leq C(K)\left( \left|
u^l-v^l\right| +\left| p-q\right| \right)
\end{array}$

\smallskip\ \\for every $(x)\in \Omega$, $\left| u\right| +\left| v\right|
+\left| p\right| +\left| q\right| \leq K$, $l=1,...N.$

Hereafter by $f^-(x)=min(f(x),0)$ and $f^+(x)=max(f(x),0)$ are
denoted the non-negative and, respectively, the non-positive part of
the function f. The same convention is valid for matrixes as well.
For instance, we denote by $M^+$ the non-negative part of $M$,
i.e.$M^+={\{ m^+_{ij}(x)\} }_{i,j=1}^{N} $.

\smallskip\ \\This paper concerns the validity of the comparison principle for
weakly-coupled elliptic systems. Let us briefly recall the
definition of the comparison principle in a weak sense for linear
systems.

{\it \ The comparison principle holds in a weak sense for the
operator $L_M$ if $(L_Mu,v)\leq 0$ and $u|_{\partial \Omega}\leq 0$
imply $(u,v)\leq 0$ in $\Omega$ for every $v>0$, $v\in
\left(W^{1,\infty}(\Omega)\cap C_0(\overline{\Omega})\right)^N $ and  $u\in
\left(W^{1,\infty}(\Omega)\cap C(\overline{\Omega})\right)^N $.}

\smallskip\ \\As it is well-known, there is no comparison principle for an
arbitrary elliptic system /see Theorem 6 below/. On the other hand,
there are broad classes of elliptic systems, such that the
comparison principle holds for their members. According to Theorem 1
below, one of these classes can be constructed using the following
condition:

\smallskip\ \\(6) {\it \qquad There is real-valued principal eigenvalue $\lambda_{\Omega_0}$
of $L_M$ and its adjoint operator ${L^*}_M$ for every
$\Omega_0\subseteq\Omega$, such that the corresponding
eigenfunctions $\tilde{w_{\Omega_0}},w_{\Omega_0}\in {\left(
W_{loc}^{2}(\Omega_0)\bigcap C_0(\overline{\Omega_0}) \right)}^N$
are positive ones}.$\Box$

\smallskip\ \\\emph{{Remark 1: By adjoint operator we mean ${L^*}_M=L^{*}+M^{t}$,
$L^{*}=diag\left(L^{*}_1, L^{*}_2,..., L^{*}_N\right) $, and
$L^{*}_{k}$ are $L^2$-adjoint operators to $L_{k}$}. The principal
eigenvalue is the first one, or the smallest eigenvalue.}

\smallskip\ \\More precisely, the class is $C^6=\{L_M$ satisfies (6) and $\lambda_{\Omega_0}>0$
for every $\Omega_0\subseteq\Omega\}$ i.e. $C^6$ contains the
elliptic systems possessing a positive principal eigenvalue with
positive corresponding eigenfunction in $\Omega_0$. In this case the
necessary and sufficient condition for the validity of the
comparison principle for systems (Theorem 1) is the same as the one
for a single equation (See [2]).

{\bf Theorem 1}{\it : Assume that (2) and (6) are satisfied. The
comparison principle holds for system (1) if the principal
eigenvalue $\lambda_{\Omega_0} >0$, where $\lambda_{\Omega_0}$ is
the principal eigenvalue of the operator $L_{M}$ on
$\Omega_0\subseteq\Omega$. If the principal eigenvalue
$\lambda=\lambda_{\Omega}\leq 0$, then the comparison principle does
not hold.}

If we consider classical solutions, then comparison principle holds
if and only if $\lambda=\lambda_{\Omega}\leq 0$.

Proof: 1.Assume that the comparison principle does not hold for
$L_M$. Let $\underline{u}, \overline{u}\in
\left(W^{1,\infty}(\Omega)\cap C(\overline{\Omega})\right)^N$ be an
arbitrary weak sub- and super-solution of $L_M$. Then
$u=\underline{u}- \overline{u}\in \left(W^{1,\infty}(\Omega)\cap
C(\overline{\Omega})\right)^N$ is a weak sub-solution of $L_M$, i.e.
$(L_M(u),v)\leq 0$ in $\Omega$ for any $v\in \left(
W^{1,\infty}(\Omega)\cap C_0(\overline{\Omega})\right)^N, v>0$ and
$u^+\equiv 0$ on $\partial\Omega$. Suppose $u^+\neq 0$. Then

$$0\geq \left(L_Mu^+,w_{\Omega_0}\right)=\left(u^+,L^*_Mw_{\Omega_0}\right)=\lambda\left(u^+,w_{\Omega_0}\right)>0$$

\smallskip\ \\for $\lambda_{\Omega_0}$, $w_{\Omega_0}$ defined in (6).

Therefore $u^+\equiv 0$, i.e for any sub- and
super-solution of $L_M$ we obtain $\underline{u}\leq \overline{u}$.

2. Suppose $\lambda\leq 0$ and $\tilde{w}$ is the corresponding
positive eigenfunction of $L_M$. Then $\tilde{w}>0$ but $L_M(\tilde{w})=\lambda
\tilde{w} \leq 0$. Therefore the comparison principle does not hold for (1).$\Box$

Unfortunately, there are some odds in the application of this
general theorem since the condition (6) is uneasy to check. First of
all, the system (1) may have no principal eigenvalue at all (See
[10]). Another obstacle is the computation of $\lambda$ even when it
exists.

Comparison principle holds for members of another broad class,
so-called cooperative elliptic systems, i.e. the systems with
$m_{ij}(x)\leq 0$ for $i\neq j$ (See [9]). Most results on the
positivity of the classical solutions of linear elliptic systems
with non-negative boundary data are obtained for the cooperative
systems (See [6,7,13,15,16,18,19,21]). As it is well known, the
positiveness and the comparison principle are equivalent for linear
systems. As for the non-linear ones, the positiveness of the
solutions is a weaker statement than the comparison principle;
positiveness can hold without ordering of sub-and super-solutions or
uniqueness of the solutions at all.

Comparison principle for the diffraction problem for weakly coupled
quasi-linear elliptic systems is proved in [3].

The spectrum properties of the cooperative $L_M$ are studied as
well. A powerful tool in the cooperative case is the theory of the
positive operators (See [17]) since the inverse operator of the
cooperative $L_{M^-}$ is positive in the weak sense. Unfortunately,
this approach cannot be applied to the general case $M\neq M^-$
since $(L_M)^{-1}$ is not a positive operator at all. Nevertheless
in [20] is proved the validity of the comparison principle for
non-cooperative systems obtained by small perturbations of
cooperative ones.

Using unconventional approach, an interesting result is obtained in
[14] for two-dimensional system (1) with $m_{11}=m_{22}=0$ and
$m_{ij}=p_i(x)>0$ for $i\neq j$, $i=1,2$. Theorem 6.5 [14] states
the existence of a principal eigenvalue with positive principal
eigenfunction in the cone $C_U=P_U\times (-P_U)$, where $P_U$ is the
cone of the positive functions in $W^1_{\infty}(\Omega)$. In the
same paper, Theorem 6.3, are provided sharp conditions for the
validity of the comparison principle with respect to the order in
$C_U=P_U\times (-P_U)$, i.e. $(u_1,u_2)\leq (v_1,v_2)$ if and only
if $u_1\leq v_1$ and $u_2\geq v_2$.

In [12] are studied existence and local stability of positive
solutions of systems with $L_k=-d_k \Delta$, linear cooperative and
non-linear competitive part, and Neumann boundary conditions.
Theorem 2.4 in [12] is similar to Theorem 2 in the present article
for $L_k=-d_k \Delta$.

Let us recall that the comparison principle was proved in [11] for
the viscosity sub-and super-solutions of general fully non-linear
elliptic systems $ G^l(x,u^1,...u^N,Du^l,D^2u^l)=0$, $l=1,...N$ /See
also the references there/. The systems considered in [11] are
degenerate elliptic ones and satisfy the same structure-smoothness
condition as the one for a single equation. The first main
assumption in [11] guarantees the quasi-monotonicity of the system.
Quasi-monotonicity in the non-linear case is an equivalent condition
to the cooperativeness in the linear one.

The second main assumption in [11] comes from the method of doubling
of the variables in the proof.

\smallskip\ \\This work extends the results obtained for cooperative systems to
the non-cooperative ones. The general idea is the separation of the
cooperative and competitive part of system (1). Then using the
appropriate spectral properties of the cooperative part, in Theorems
3 and 4 are derived conditions for the validity of the comparison
principle for the initial system. In particular in Theorem 3 is
employed the fact that irreducible cooperative system possesses a
principal eigenvalue and the corresponding eigenfunction is a
positive one, i.e. condition (6) holds. This way are obtained some
sufficient conditions for validity of the comparison principle for
the non-cooperative system as well. Analogously, in Theorem 4 are
derived the corresponding conditions for the validity of comparison
principle for competitive systems. The conditions derived in
Theorems 3 and 4 are not sharp.

Since predator-prey systems are basic model example for
non-cooperative systems, in Theorem 5 is adapted the main idea of
Theorem 4 to systems which cooperative part is a triangular matrix.
Sufficient condition for the validity of comparison principle for
predator-prey systems is derived in Theorem 5.

In Theorems 6 and 7 are given conditions for failure of the
comparison principle.

The results of Theorems 3 and 4 are adapted to quasi-linear systems
in Theorem 8.

\section{Comparison principle for linear elliptic systems}

As a preliminary statement we need the following well known fact

{\bf Theorem 2}{\it : Every irreducible cooperative system $L_{M^-}$
has unique principal eigenvalue and the corresponding eigenfunction
is positive }.

The principal eigenfunction for linear operators is unique up to
positive multiplicative constants, but for our purpose the
positiveness is of importance.

In fact, Theorem 2 is in the scope of Theorems 11 and 12 in [1].
Theorems 11 and 12 in [1] concern second order cooperative linear
elliptic systems with cooperative boundary conditions and are more
general then Theorem 2. In sake of completeness, a sketch of the
proof of Theorem 2 follows. It is based on the idea of adding a big
positive constant to the operator. The same idea appears for
instance in [16] and many other works.

Skatch of tne proof: Let us consider the operator $L_c= L_{M^-}+cI$
where $c\in R$ is a constant and $I$ is the identity matrix in
$R^n$. Then $L_c$ satisfies the conditions of Theorem 1.1.1 [16] if
$c$ is large enough, namely

1.  $L_c$ is a cooperative one;

2. $L_c$ is a fully coupled;

3. There is a super-solution $\varphi$ of $L_c \varphi =0$.

Conditions 1 and 2 above are obviously fulfilled by $L_c$, since
$L_{M^-}$ is a cooperative and a fully coupled one, and $L_c$
inherits these properties from $L_{M^-}$.

As for the condition 3, we construct the super – solution $\varphi$
using the principal eigenfunctions of the operators $L_k-c_k$. More
precisely, $\varphi = (\varphi_1, \varphi_2,..., \varphi_N)$, where
$ \left( L_k-c_k\right) \varphi_k={\lambda}_k \varphi_k$, and
${\lambda}_k,\varphi_k>0$ in $\Omega$. The existence of $\varphi_k$
is a well - known fact.

We claim that if $c$ is large enough then $\varphi$ is a super -
solution
 of $L_c$ , i.e. $\varphi \in {\left( W_{loc}^{2,n}(\Omega)\bigcap
C(\overline{\Omega})\right)}^N$ and $\varphi\geq 0$, $L_c \varphi
\geq 0$ and $\varphi$ is not identical to null in $\Omega$.

Since we have chosen ${\varphi}_k$ being the principal
eigenfunctions of $L_k-c_k$, we have ${\varphi}_k \in {\left(
C^{2}(\Omega)\bigcap C(\overline{\Omega})\right)}$ and ${\varphi}_k>
0$. It remains to prove that $L_c \varphi \geq 0$.

Let $$A_k = {\left(L_c\varphi \right)}_k = -\sum_{i,j=1}^{n}D_j
\left( a_k^{ij}(x)D_i{\varphi}_k \right)
+\sum_{i=1}^{n}b_k^i(x)D_i{\varphi}_k+
\sum_{i=1}^{n}m_{ki}(x){\varphi}_i +(c_k+c){\varphi}_k =$$

$$=({\lambda}_k+c_k+c){\varphi}_k+
\sum_{i=1}^{n}m_{ki}(x){\varphi}_i.
$$

Then $A_k\geq 0$ for every $i$.

First of all, if we denote by $n$ the outer unitary normal vector to
$\partial\Omega$, then
$${\frac{dA_k}{dn}}|_{\partial\Omega}=
({\lambda}_k+c_k+c)\frac{d{\varphi}_k}{dn}
+\sum_{i=1}^{n}m_{ki}(x)\frac{d{\varphi}_i}{dn}$$ since
${\varphi_i}|_{\partial\Omega}=0$. Therefore there is a constant
$c'$, such that ${\frac{dA_k}{dn}}|_{\partial\Omega}<0$ for $c>c'$
since $\frac{{d\varphi}_i}{dn}<0$ on $\partial\Omega$ (See [14],
Theorem 7, p.65) and $\lambda_i$ is independent on $c$.

Hence there is a neighbourhood ${\Omega}_\varepsilon =
\{x\in\overline{\Omega}:dist(x,\partial\Omega)<\varepsilon\}$ for
some $\varepsilon>0$, such that
$${\frac{dA_k}{dn}}|_{{\Omega}_\varepsilon}<0$$.

Since $A_k=0$ on $\partial\Omega$, then $A_k>0$ in
${\Omega}_\varepsilon $

The set $\Omega \setminus {\Omega }_\varepsilon$ is compact,
therefore there is $c">0$ such that $A_k>0$ in the compact set
$\Omega \setminus {\Omega}_\varepsilon$ for $c>c"$, since
${\varphi}_k>0$ in ${\Omega} \setminus {\Omega}_{\varepsilon}$.

Considering $c>max(c',c")$ we obtain $A_k>0$ in $\Omega$, therefore
$\varphi$ is indeed a super - solution of $L_c$.

The rest of the proof follows the proof of Theorem 1.1.1 [16].$\Box$

A reasonable question is: could the non-cooperative part of the
system "improve" the spectral facilities of the cooperative system?
In other words, if the cooperative part of the system has
non-positive principal eigenvalue, what are conditions on the
competitive part, such that the comparison principle holds for the
system? An answer of this question is given in the following

{\bf Theorem 3}{\it : Let (1) be a weakly coupled system with
irreducible cooperative part of $L^*_{M^-}$ such that (2) is
satisfied. Then the comparison principle holds for system (1) if
there is $x_0\in \Omega$ such that

\smallskip\ \\(7)\qquad $%
 \left(\lambda+\sum_{k=1}^{N}m_{kj}^{+}(x_0)\right)>0$ for
 $j=1...N$

\smallskip\ \\and

\smallskip\ \\(8)\qquad $%
 \lambda+ m_{jj}^{+}(x)\geq 0$ for every $x\in \Omega$ and
 $j=1...N$

\smallskip\ \\where $\lambda=inf_{\Omega_0\subseteq \Omega}\{\lambda_{\Omega_0}$ : $\lambda_{\Omega_0}$ is the
principal eigenvalue of the operator $L_{M^-}$ on $\Omega_0$\}}.

It is obvious, that if $ \lambda> 0$, then the comparison principle
holds. More interesting case is $ \lambda < 0$. Then $m^+_{kj}$ can
"improve" the properties of $L_M$ with respect to the validity of
the comparison principle. Furthermore, if $ \lambda+ m_{jj}^{+}(x)>
0$, then (7) is consequence of (8). Condition (7) is important when
$ \lambda+ m_{jj}^{+}(x)\equiv 0$.

\emph{Remark 2: If $L^*_{M^-}$ is irreducible, then $L_{M^-}$ is
irreducible as well. In fact $L^*_{M^-}=L^*+{M^-}^t$ and if
${M^-}^t$ is irreducible, then such is $M^-$.}

Proof: Suppose all conditions of Theorem 3 are satisfied by $L_M$
but the comparison principle does not hold for $L_M$. Let
$\underline{u}, \overline{u}\in \left(W^{1,\infty}(\Omega)\cap
C(\overline{\Omega})\right)^N$ be an arbitrary weak sub- and
super-solution of $L_M$. Then $u=\underline{u}- \overline{u}\in
\left(W^{1,\infty}(\Omega)\cap C(\overline{\Omega})\right)^N$ is a
weak sub-solution of $L_M$ as well, i.e. $(L_M(u),v)\leq 0$ in
$\Omega$ for any $v\in \left(W^{1,\infty}(\Omega)\cap
C_0(\overline{\Omega}\right)^N, v>0$ and $u^+\equiv 0$ on
$\partial\Omega$.

Assume $u^+\neq 0$. Let $\Omega_{supp(u^+)}\subseteq supp(u^+)$ has
smooth boundary. Then for any $v> 0$, $v\in
\left(W^{1,\infty}(\Omega_{supp(u^+)})\cap
C(\overline{\Omega_{supp(u^+)}})\right)^N$

\smallskip\ \\(9)\qquad $0\geq \left( \L_{M}u^+,v\right)
=\left(u^+,\L^*_{M^-}v\right)+\left(M^+u^+,v\right)$

\smallskip\ \\is satisfied since $L_{M}(u^+)\leq 0$.

Since $L_{M^-}$ is a cooperative operator, such is ${\left(
L_{M^-}\right)}^{*}=L^{*}+(M^-)^{t}$ as well. According to Theorem 2
above, there is a unique positive eigenfunction $w\in {\left(
W_{loc}^{2,n}(\Omega_{supp(u^+)})\bigcap
C_0(\overline{\Omega_{supp(u^+)}}) \right)}^N$ such that $w>0$ and
$L^*_{M^-}w=\lambda w$ for some $\lambda >0$.

Then $w$ is a suitable test-function for (9). Rewriting the
inequality (9) for $v=w$ we obtain

\smallskip\ \\$0\geq \left(u^+,\L^*_{M^-}w\right)+\left(M^+u^+,w\right)=
\left(u^+,\lambda w\right)+\left(M^+u^+,w\right)$

\smallskip\ \\or componentwise

\smallskip\ \\(10)\qquad
$0\geq \left(u^{+}_{k},\lambda
w_k\right)+\left({\sum}_{j=1}^{N}m^{+}_{kj}u^{+}_{j},w_k\right)$

\smallskip\ \\for $k=1,...n$.

The sum of inequalities (10) is

\smallskip\ $0\geq {\sum}_{k=1}^{N}\left(
\left(u^{+}_{k},{{\tilde{L}}^*}_{k}w_k\right)
+\left({\sum}_{j=1}^{N}m^{+}_{kj}u^{+}_{j},w_k\right)\right)=$

\smallskip\ \\
$={\sum}_{k=1}^{N}\left(u^{+}_{k},{\lambda}w_k\right)
+{\sum}_{k,j=1}^{N}\left(u^{+}_{j},m^{+}_{kj}w_k\right)=$

\smallskip\ \\
$={\sum}_{j=1}^{N}\left(u^{+}_{j},{\sum}_{k=1}^{N}\left(
{{\delta}_{jk}\lambda}+m^{+}_{kj}\right)w_k\right)>0$

\smallskip\ \\since $u^{+}> 0$, $w_k>0$, (7) and (8). Condition (8)
is used in $\left(u^{+}_{k},({\lambda}+m^{+}_{kk})w_k\right)\geq 0$.

\smallskip\ \\The above contradiction proves that $u^+\equiv 0$
 and therefore the comparison principle holds for operator
$L_M$.$\Box$

Since in [1] and [18] are considered only systems with irreducible
cooperative part, the ones with reducible $L_{M^-}$ are excluded of
the range of Theorem 3. Nevertheless the same idea is applicable to
some systems with reducible cooperative part as well, as it is given
it Theorem 4.

{\bf Theorem 4}{\it : Assume $m^-_{ij}\equiv 0$ for $i\neq j$ and
(2) is satisfied. Then the comparison principle holds for system (1)
if there is $x_0\in \Omega$ such that

\smallskip\ \\(11)\qquad $%
 \left(\lambda_{j}+\sum_{k=1}^{N} m_{kj}^{+}(x_0)\right)>0$ for
 $j=1...N$

\smallskip\ \\and

\smallskip\ \\(12)\qquad $%
 \lambda_{j}+ m_{jj}^{+}(x)\geq 0$ for every $x\in \Omega$
 $j=1...N$,

\smallskip\ \\where $\lambda_{j}=inf_{\Omega_0\subseteq\Omega}\{\lambda_{j\Omega_0}$ : $\lambda_{j\Omega_0}$ is the
principal eigenvalue of the operator $L_j+m^-_{jj}$ on
$\Omega_0$\}}.

Theorem 4 is formulated for diagonal matrix $M^-$. The statement is
valid with obvious modification if $M^-$ has block structure, i.e.

$$M^-=\left(
\begin {array}{cccccccc}
M^-_1 &  & 0 &  & ... &  & 0  \\
&  &  &  &  &  &   \\
 0 &  & M^-_2 &  & ... &  & 0
 \\  &  &  &  &  &  &
\\... &  & ... &  & ... &  & ... & \\
&  &  &  &  &  & \\
0 &  & 0 &  & ... &  &M^-_r
\end{array}
\right) $$

\smallskip\ \\where $M^-_k$ are $d_k$-dimensional square matrixes,
$\sum d_k\leq N$ .

Proof: Let all conditions of Theorem 4 be satisfied by $L_M$ but the
comparison principle does not hold for $\tilde{L}_{M^+}$. Let
$\underline{u}, \overline{u}\in \left(W^{1,\infty}(\Omega)\cap
C(\overline{\Omega})\right)^N$ be an arbitrary weak sub- and
super-solution of $\tilde{L}_{M^+}$. Then $u=\underline{u}-
\overline{u}\in \left(W^{1,\infty}(\Omega)\cap
C(\overline{\Omega})\right)^N$ is a weak sub-solution of
$\tilde{L}_{M^+}$ as well, i.e. $(\tilde{L}_{M^+}(u),v)\leq 0$ in
$\Omega$ for any $v\in \left(W^{1,\infty}(\Omega)\cap
C_0(\overline{\Omega})\right)^N, v>0$ and $u^+\equiv 0$ on
$\partial\Omega$.

Suppose that $u^+\neq 0$. Let $\Omega_{supp(u^+)}\subseteq
supp(u^+)$ has smooth boundary. Then for any $v> 0$, $v\in
W_2^{1,\infty}(\Omega_{supp(u^+)})\cap
C(\overline{\Omega_{supp(u^+)}})$

\smallskip\ \\(13)\qquad $0\geq \left( \tilde{L}_{M^+}u^+,v\right)
=\left(u^+,\tilde{L}^*v\right)+\left(M^+u^+,v\right)$

\smallskip\ \\is satisfied since $\tilde{L}_{M^+}u^+\leq 0$.

According to Theorem 2.1 in [2], there is a positive principal
eigenfunction for the operator ${\tilde{L}}^*_k$ in
$\Omega_{supp(u^+)}$, i.e. $\exists\qquad w_k(x)\in
C^2(\Omega_{supp(u^+)} \bigcap R^1)$ such that
${{\tilde{L}}^*}_kw_k(x)={\lambda}_kw_k(x)$ and $w_k(x)>0$. Note
that $w_k$ are classical solutions.

Then the vector-function $w(x)=(w_1(x),...,w_n(x))$, composed of the
principal eigenfunctions $w_k(x)$, is suitable as a test-function in
(13).

Writing componentwise inequality (13) for $v=w$ we obtain

\smallskip\ \\(14)\qquad
$0\geq
\left(u^{+}_{k},{{\tilde{L}}^*}_{k}w_k\right)+\left({\sum}_{j=1}^{N}m^{+}_{kj}u^{+}_{j},w_k\right)$

\smallskip\ \\for $k=1,...N$.

The sum of inequalities (14) is

\smallskip\ $0\geq {\sum}_{k=1}^{N}\left(
\left(u^{+}_{k},{{\tilde{L}}^*}_{k}w_k\right)
+\left({\sum}_{j=1}^{N}m^{+}_{kj}u^{+}_{j},w_k\right)\right)=$

\smallskip\ \\
$={\sum}_{k=1}^{N}\left(u^{+}_{k},{\lambda}_{k}w_k\right)
+{\sum}_{k,j=1}^{N}\left(u^{+}_{j},m^{+}_{kj}w_k\right)=$

\smallskip\ \\
$={\sum}_{j=1}^{N}\left(u^{+}_{j},{\sum}_{k=1}^{N}\left(
{{\delta}_{jk}\lambda}_{j}+m^{+}_{kj}\right)w_k\right)>0$

\smallskip\ \\since $u^{+}> 0$, $w_k>0$, (11) and (12).

\smallskip\ \\The above contradiction proves that $u^+\equiv 0$
 and therefore the comparison principle holds for operator
$L_M^{+}$.$\Box$

\emph{Remark 3: It is obvious that conditions (7),(8), and
respectively, (11), (12), can be substituted by the sharper
condition ${\sum}_{k=1}^{n}\left(
{{\delta}_{jk}\lambda_k}+m_{kj}\right) w_k>0$ for every $x\in
\Omega$ and every $j=1...N$, which is useful only if the exact
values of the eigenfunctions $w_k$ can be computed.}

The main idea in Theorem 4 could be modified for systems with
triangular cooperative part, for instance with null elements above
the main diagonal. For instance predator-prey systems have
triangular cooperative part. Of course, if $m^-_{ij}(x)>0$ for every
$x\in\Omega$ and $i=1,...N$, $j<i$, then the system is in the scope
of Theorem 3. In Theorem 5 this condition is not necessary, i.e.
some of the species can extinguish in some subarea of $\Omega$.

{\bf Theorem 5}{\it : Assume (2) is satisfied and the cooperative
part $M^-$ is triangular for the system (1), i.e. $m^-_{ij}=0$ for
$i=1,...N$, $j>i$. Then the comparison principle holds for system
(1), if there is} $\varepsilon > 0$ \emph{such that}

\smallskip\ \\(15)\qquad $%
 \left(\lambda_{j}-(1-\delta_{1j})\varepsilon+\sum_{k=1}^{N} m_{kj}^{+}(x_0)\right)>0$ for
 $j=1...N$ \emph{for some} $x_0\in\Omega$

\smallskip\ \\\emph{and}

\smallskip\ \\(16)\qquad $%
 \lambda_{j}-(1-\delta_{1j})\varepsilon +m^+_{jj}(x)\geq 0$ for
 every $x\in\Omega$ and $j=1...N$,

\smallskip\ \\\emph{where $\lambda_{j}=inf_{\Omega_0\subseteq\Omega}\{\lambda_{j\Omega_0}$ : $\lambda_{j\Omega_0}$ is the
principal eigenvalue of the operator $L_j+m^-_{jj}$ on
$\Omega_0$\}}.

Note that the condition for triangular cooperative part  does not
exclude $m^-_{ij}(x_0)=0$ for some $x_0\in \Omega$, $i,j=1,...N$.

Proof: 1. The first equation in $L_{M^-}$ is not coupled, and there
are principal eigenvalue $\lambda_1$ and principal eigenfunction
$w_1>0$ of $L_1+m^-_{11}$ (See Theorem 2.1 in [2]). We put
$\widetilde{w}_1=w_1$.

2. The equation
$(L_2+m^-_{22})\widetilde{w}_2-\lambda\widetilde{w}_2 =
m_{21}\widetilde{w}_1$ with null boundary conditions has unique
solution for $\lambda<\lambda_2$, where $\lambda_2$ is the principal
eigenvalue of $L_2+m^-_{22}$. We put
$\lambda=\lambda_2-\varepsilon$. Since the right-hand side
$m_{21}\widetilde{w}_1$ is positive, the solution $\widetilde{w}_2$
is positive as well.

3. By induction we construct positive functions $\widetilde{w}_j$,
$j=3,...N$ as solutions of
$(L_j+m^-_{jj})\widetilde{w}_j-(\lambda_j-\varepsilon)\widetilde{w}_j
= \sum_{i=1}^{j-1}m_{ji}\widetilde{w}_i$ with null boundary
conditions. As usual $\lambda_j$ are the principal eigenfunctions of
$L_j+m^-_{jj}$.

4. The rest of the proof follows the proof of Theorem 4 where
$\lambda_j$ is substituted with $\lambda_j-\varepsilon$ and $w_j$ is
substituted with $\widetilde{w}_j$.

For the simplest predator-prey system, $N=2$, $m_{11}=m_{22}=0$,
$m_{12}>0$ and $m_{21}<0$, conditions (15) and (16) are
$\lambda_1\geq0$, $\lambda_2>0$, where $\lambda_{j}$ is the
principal eigenvalue of the operator $L_j$, $j=1,2$.

Condition (12) in Theorem 2 is useful for construction of
counter-example for the non-validity of comparison principle in
general.

{\bf Theorem 6}{\it :  Let (1) be a weakly coupled system with
reducible cooperative part $L_{M^-}$ and (2) be satisfied. Suppose
that (12) is not true, i.e there is some $j\in \{ 1...N\}$ such that
$ \left( \lambda_j + m_{jj}^{+}(x)\right)<0$ for any $x\in\Omega$,
and $m^+_{jl}=0$ for $l\neq j$, $l=1,...N$. Then comparison
principle does not hold for system (1)}.

Proof: Let us suppose for simplicity that $j=1$ and $m^-_{1,j}=0$
for $j=2,...N$. We consider vector-function $w(x)={w_1(x),0,...,0}$,
where $w_1(x)$ is the principal eigenfunction of  $L_1+m^-_{11}$.

Then for the first component ${(L_M)}_1$ of $L_M$ is valid
$(L_Mw)_1= \lambda w_1(x) + m_{11}^{+}w_1(x)<0$ in  $\Omega$, where
$\lambda_j$ is the principal eigenvalue of $L_1$, and $(L_Mw)_k=0$
for $k=1,...N$. Therefore, $L_Mw\leq 0$ but $w(x)\geq 0$ and
comparison principle fails. $\Box$

The simplest case to illustrate Theorems 4 and 6 is $N=2$. Let us
consider irreducible competitive system

\smallskip\ \\(17)\qquad $L_ju_j+\sum_{j,k=1}^2m_{jk}u_k=f_j$, $j=1,2$,

\smallskip\ \\where $m_{11}=m_{22}=0$, $m_{12}>0$, $m_{21}>0$.

Suppose $\lambda_ j$ is the principal eigenvalue of $L^*_j$,
$j=1,2$. If $\lambda_ j\geq 0$ and there is $x_0\in\Omega$  such
that $\lambda_1+m_{21}(x_0)>0$  and $\lambda_2+m_{12}(x_0)>0$, then
according to  Theorem 4 the comparison principle holds for system
(1), i.e. if  $f_1>0$, $f_2>0$, then $u_1>0$ and $u_2>0$, where
$u=\underline{u}-\overline{u}$ is defined in the proof of Theorem 3.

If $\lambda_2+m_{12}(x)<0$ for every $x\in\Omega$, then according to
Theorem 6 there is no comparison principle for system (1) in the
lexicographic order, used in this paper.

More detailed analysis of the validity of the comparison principle
for system (1) could be done if we consider order in the cone
$C_U=P_U\times (-P_U)$, i.e. $(u_1,u_2)\leq (v_1,v_2)$ if and only
if $u_1\leq v_1$ and $u_2\geq v_2$. Then Theorem 6.5 [14] states the
existence of a principal eigenvalue $\lambda$ of $L^*$ with positive
in $C_U$ principal eigenfunction $w_1(x)>0$, $w_2(x)<0$.

If $\lambda>0$, then according to Theorem 6.3 [14] the comparison
principle holds in the order in $C_U$, i.e. if $f_1>0$, $f_2<0$,
then $u_1>0$ and $u_2<0$.

If $\lambda<0$, then
$(L_1(-u_1)+m_{12}u_2,w_1)+(L_2u_2+m_{21}(-u_1),w_2)= (-u_1, \lambda
w_1+m_{21}w_2)+(u_2, m_{21}w_1+\lambda w_2)>0$. Hence $u_1<0$ and
$u_2>0$ for $f_1>0$, $f_2>0$.

A statement analogous to Theorem 6 is valid for irreducible systems
as well.

{\bf Theorem 7}{\it :  Let (1) be a weakly coupled system with
irreducible cooperative part $L_{M^-}$ and (2) be satisfied. Suppose
that (7) is not true, i.e there is some $j\in \{ 1...N\}$ such that
$ \left( \lambda + m_{jj}^{+}(x)\right)<0$ for any $x\in\Omega$, and
$m^+_{jl}=0$ for $l\neq j$, $l=1,...N$. Then comparison principle
does not hold for system (1)}.

Note that in Theorem 6 and Theorem 7 we need the violation of
condition (12) and, respectively,  condition (7) in all $\Omega$.
The proof of Theorem 7 follows the proof of Theorem 6 with obvious
adaptation.

\section{Comparison principle for quasi-linear elliptic systems}

Considering quasi-linear system (3), (4), we use the results of
the previous section to derive conditions for the validity of comparison principle.

Let $u(x)\in \left(W^{1,\infty}(\Omega)\cap
C(\overline{\Omega})\right)^N$ be a sub-solution and $v(x)\in
\left(W^{1,\infty}(\Omega)\cap C(\overline{\Omega})\right)^N$ be a
super-solution of (3), (4). Comparison principle holds for (3), (4),
if $Q(u)\leq Q(v)$ in $\Omega$, $u\leq v$ on $\partial \Omega$ imply
$u\leq v$ in $\Omega$. Last three inequalities are considered in the
weak sense.

Recall that the vector-function $u(x)$ is a weak sub-solution of (3), (4) if

$$\int_\Omega \left( a^{li}(x,u^l,Du^l)\eta _{x_i}^l+F^l(x,u^1,...u^N,Du^l)\eta ^l-f^l(x)\eta ^l\right) dx\leq 0$$

\smallskip\ \\for $l=1,...N$ and for every nonnegative vector function
$\eta\in \left(\stackrel{\circ }{W}^{1,\infty}(\Omega)\cap
C(\overline{\Omega})\right)^N$ (i.e. $\eta=(\eta ^1,...\eta ^N)$,
$\eta ^l\geq 0$, $\eta ^l\in \left(W^{1,\infty}(\Omega)\cap
C(\overline{\Omega})\right)^N\cap C( \overline{\Omega})$ and  $\eta
^l=0$ on $\partial \Omega$.

Analogously, $v(x)\in \left(W^{1,\infty}(\Omega)\cap
C(\overline{\Omega})\right)^N$ is a super-solution of (3), (4), if
$$\int_\Omega \left( a^{li}(x,v^l,Dv^l)\eta _{x_i}^l+F^l(x,v^1,...v^N,Dv^l)\eta ^l-f^l(x)\eta ^l\right) dx\geq 0$$

\smallskip\ \\for $l=1,...N$ and for every nonnegative vector function
$\eta \in \left(\stackrel{\circ }{W}^{1,\infty}(\Omega)\cap
C(\overline{\Omega})\right)^N$.

Since $u(x)$ and $v(x)$ are sub-and super-solution respectively,
then $\tilde{w}(x)=u(x)-v(x)$ is a weak sub-solution of the
following problem

\smallskip\ \\ $-\sum^{n}_{i,j=1} D_i \left( B_j^{li}D_j \tilde{w}^l%
+B_0^{li}\tilde{w}^l\right) +{\sum}^{N}_{k=1}E_k^l\tilde{w}^k+{\sum}^{n}_{i=1}H_i^lD_i \tilde{w}^l=0$ in $%
\Omega $

\smallskip\ \\with non-positive boundary data on $\partial \Omega $. Here

\smallskip\ \\$B_j^{li}=\int_0^1\frac{\partial a^{li}}{\partial p_j}%
(x,P^l)ds$, $B_0^{li}=\int_0^1\frac{\partial a^{li}}{\partial
u^l}(x,P^l)ds$, $E_k^l=\int_0^1\frac{\partial F^l}{\partial
u^k}(x,S^l)ds$,

\smallskip\ \\$H_i^l=\int_0^1\frac{\partial F^l}{\partial p_i}(x,S^l)ds$, $%
P^l=\left( v^l+s(u^l-v^l),Dv^l+sD(u^l-v^l)\right) $,

\smallskip\ \\ $S^l=\left( v+s(u-v),Dv^l+sD(u^l-v^l)\right)$.

Therefore, $\tilde{w}_{+}(x)=\max \left( \tilde{w}(x),0\right)$ is a
sub-solution of

\smallskip\ \\(18)\qquad $\sum^{n}_{i,j=1} D_i\left( B_j^{li}D_j{\tilde{w}_+}^l%
+B_0^{li}{\tilde{w}_+}^l\right) +\sum^{N}_{k=1}E_k^l{\tilde{w}_+}^k+\sum^{n}_{i=1}H_i^lD_i {\tilde{w}_+}^l=0$ in $%
\Omega $

\smallskip\ \\with zero boundary data on $\partial \Omega$.

Equation (18) is equivalent in terms of matrix to

\smallskip\ \\(19)\qquad $B_E\tilde{w}_+=(B+E)\tilde{w}_+=0$ in
$\Omega$,

\smallskip\ \\where $B=diag(B_1,B_2,...B_N)$,
$B_l=\sum^{n}_{i,j=1} D_i\left( B_j^{li}D_j{\tilde{w}_+}^l%
+B_0^{li}{\tilde{w}_+}^l\right) +\sum^{n}_{i=1}H_i^lD_i
{\tilde{w}_+}^l$ and $E=\{E_k^l\}_{l,k+1}^N$.

If we denote $B_i^{kj}$ by $a_k^{ij}$, $B_0^{ki}+H_i^k$ by $b_k^i$,
${\sum}_{i=1}^n D_iB_0^{ki}+E_k^k$ by $m_{kk}(x)$ for $i,j=1...n$,
$k=1...N$ and $E_k^l$ by $m_{lk}(x)$ for $k,l=1...N$, $k\neq l$,
system (18) looks like system (1). Hereafter we follow the notations
for system (1).

Suppose now that $\tilde{w}_{+}(x)$ is not identical equal to zero
in $\Omega $, i.e. comparison principle fails for (3), (4). Suppose
$L_{M^-}$ is irreducible. Then

\smallskip\ \\$0\geq \left(L_M \tilde{w}_+,w \right)
=\left(\tilde{w}_{+},\L^*_{M^-}w\right)+\left(M^+\tilde{w}_{+},w\right)=
\left(\tilde{w}_{+},\lambda
w\right)+\left(M^+\tilde{w}_{+},w\right)$

\smallskip\ \\where $\lambda$ is the principal eigenvalue of $L^*_{M^-}$ and $w$ is the corresponding eigenfunction.

Suppose $a_k^{ij}$ and $m_{lk}(x)$ satisfy the conditions (2), (7)
and (8) in Theorem 3. Following the proof of Theorem 3, we obtain
that $\tilde{w}_{+}\equiv 0$ in $\Omega$, i.e. comparison principle
holds for the system (3), (4).

If $L_{M^-}$ is reducible, then

\smallskip\ \\$0\geq \left(L_M\tilde{w}_{+},w \right)
=\left(\tilde{w}_{+},\L^*w\right)+\left(M^+\tilde{w}_{+},w\right)=\left(\tilde{w}_{+},\tilde{\lambda}
w\right)+\left(M^+\tilde{w}_{+},w\right)$

\smallskip\ \\where
$\tilde{\lambda}w=(\tilde{\lambda}_1w_1,\tilde{\lambda}_2w_2,...\tilde{\lambda}_Nw_N)$,
$\tilde{\lambda}_k$ is the principal eigenvalue of $L^*_{k}$ and
$w_k$ is the corresponding eigenfunction for $k=1,...N$.

Suppose $a_k^{ij}$ and $m_{lk}(x)$ satisfy the conditions (2), (11)
and (12) in Theorem 4. Following the proof of Theorem 4, we obtain
that $\tilde{w}_{+}\equiv 0$ in $\Omega$, i.e. comparison principle
holds for the system (3), (4).

We have sketched the proof the following

{\bf Theorem 8}{\it :  Suppose (3), (4) is a quasi-linear system and
the corresponding system $B_{E^-}$ in (19) is elliptic. Then the
comparison principle holds for system (3), (4) if}

\smallskip\ \\(i)\qquad {\it $B_{E^-}$ in (19) is irreducible and for every $j=1...n$}

\smallskip\ \\(ii)\qquad $%
\lambda+\left(\sum_{k=1}^{N} \frac{\partial F^k}{\partial
p^j}(x,p,Dp^l)+{\sum}_{i=1}^N D_i\frac{\partial a^{ji}}{\partial
p^j}(x,p^j,Dp^j)\right)^+ >0$,

\smallskip\ \\(iii)\qquad $%
 \lambda+ \left({\sum}_{i=1}^n D_i\frac{\partial a^{ji}}{\partial
p^j}(x,p^j,Dp^j)+{\frac{\partial F^j}{\partial
p^j}(x,p,Dp^j)}\right)^+\geq 0$

\smallskip\ \\ {\it where $x\in \Omega$, $p\in R^n$ and $\lambda=inf_{\Omega_0\subseteq \Omega}\{\lambda_{\Omega_0}$ : $\lambda_{\Omega_0}$ is the
principal eigenvalue of the operator $B_{E^-}$ on $\Omega_0$\};

\smallskip\ \\{\it or

\smallskip\ \\(i')\qquad $B_{E^-}$ in (19) is reducible and for every $j=1...n$}

\smallskip\ \\(ii')\qquad $%
 {\lambda}_j+\left(\sum_{k=1}^{n} \frac{\partial F^k}{\partial
p^j}(x,p,Dp^j)+{\sum}_{i=1}^n D_i\frac{\partial a^{ji}}{\partial
p^j}(x,p^j,Dp^j)\right)^+ >0$,

\smallskip\ \\(iii')\qquad $%
 {\lambda}_j+ \left({\sum}_{i=1}^n D_i\frac{\partial a^{ji}}{\partial
p^j}(x,p^j,Dp^j)+{\frac{\partial F^j}{\partial
p^j}(x,p,Dp^j)}\right)^+\geq 0$

\smallskip\ \\{\it where $x\in \Omega$, $p\in R^n$ and $\lambda_{l}=inf_{\Omega_0\subseteq\Omega}\{\lambda_{l\Omega_0}$ : $\lambda_{l\Omega_0}$ is the
principal eigenvalue of the operator $B_l$ on $\Omega_0$.}

\section{Final remarks}

The sufficient conditions in Theorems 3 and 4 are derived from the
spectral properties of the cooperative part of (1) - the operator
$L_{M^-}$, or, in other words, comparing the principal eigenvalue of
$L_{M^+}$ with the quantities in $M^{+}$. In fact the positive
matrix $M^+$ causes a migration of the principal eigenvalue of
$L_{M^-}$ to the left.

Theorems 3 and 4 provide a huge class of non-cooperative systems
such that the comparison principle is valid for. The idea of
migrating the spectrum of a positive operator on the right works in
this case, though the spectrum itself is not studied in this
article. The results for non-cooperative systems in this paper are
not sharp and the validity of the comparison principle is to be
determined more precisely in the future.

\section{Acknowledgment}

The author would like to acknowledge Professor Alexander Sobolev for
the very useful talks on the theory of positive operators, during
the author's stay at University of Sussex as Maria Curie fellow.

\section{REFERENCES}

[1] H.Amann, Maximum Principles and Principal Eigenvalues, 10
Mathematical Essays on Approximation in Analysis and Topology
(J.Ferrera, J.Lopez-Gomez and F.R.Ruiz del Portal Eds.), Elsevier,
Amsterdam (2005), 1-60.

  [2] H.Berestycki, L.Nirenberg, S.R.S. Varadhan : The principal
eigenvalue and maximum principle for second-order elliptic operators
in general domains, Commun. Pure Appl. Math. 47, No.1, 47-92 (1994).

[3] G.Boyadzhiev, N.Kutev : Diffraction problems for quasilinear
reaction-diffusion systems, Nonlinear Analysis 55 (2003), 905-926.

[4] G.Caristi, E. Mitidieri : Further results on maximum principle
for non-cooperative elliptic systems. Nonl.Anal.T.M.A., 17 (1991),
547-228.

[5] C.Coosner, P.Schaefer : Sign-definite solutions in some linear
elliptic systems. Peoc.Roy.Soc.Edinb.,Sect.A 111, (1989), 347-358.

[6] D.di Figueredo, E.Mitidieri : Maximum principles for cooperative
elliptic systems. C.R.Acad.Sci. Paris, Ser. I, 310 (1990), 49-52.

[7] D.di Figueiredo, E.Mitidieri : A maximum principle for an
elliptic system and applications to semi-linear problems, SIAM
J.Math.Anal. 17 (1986), 836-849.

[8] Gilbarg, D and Trudinger, N. Elliptic partial differential
equations of second order. 2nd ed., Springer - Verlag, New York.

[9] M.Hirsch : Systems of differential equations which are
competitive or cooperative I. Limit sets, SIAM J. Math. Anal. 13
(1982), 167-179.

[10] P.Hess : On the Eigenvalue Problem for Weakly Coupled Elliptic
Systems, Arch. Ration. Mech. Anal. 81 (1983), 151-159.

[11] Ishii, Sh. Koike : Viscosity solutions for monotone systems of
second order elliptic PDEs. Commun. Part.Diff.Eq. 16 (1991), 1095 -
1128.

[12] Li Jun Hei, Juan Hua Wu : Existence and Stability of Positive
Solutions
 for an Elliptic Cooperative System. Acta Math. Sinica Oct.2005, Vol.21,
No 5, pp 1113-1130.

[13] J.Lopez-Gomez, M. Molina-Meyer : The maximum principle for
cooperative weakly coupled elliptic systems and some applications.
Diff.Int.Eq. 7 (1994), 383-398.

[14]  J.Lopez-Gomez, J.C.Sabina de Lis, Coexistence states and
global attractivity for some convective diffusive competing species
models, Trans.Amer.Math.So. 347, 10 (1995), 3797-3833.

[15] E.Mitidieri, G.Sweers : Weakly coupled elliptic systems and
positivity. Math.Nachr. 173 (1995), 259-286.

[16] M. Protter, H.Weinberger : Maximum Principle in Differential
Equations, Prentice Hall, 1976.

[17] M.Reed, B.Simon : Methods of modern mathematical Physics, v.IV:
Analysis of operators, Academic Press, New York, (1978).

[18] G.Sweers : Strong positivity in $C(\overline{\Omega })$ for
elliptic systems. Math.Z. 209 (1992), 251-271.

[19] G.Sweers : Positivity for a strongly coupled elliptic systems
by Green function estimates. J Geometric Analysis, 4, (1994),
121-142.

[20] G.Sweers : A strong maximum principle for a noncooperative
elliptic systems. SIAM J. Math. Anal., 20 (1989), 367-371.

[21] W.Walter : The minimum principle for elliptic systems.
Appl.Anal.47 (1992), 1-6.

Author's address:

Institute of Mathematics and Informatics,

Bulgarian Academy of Sciences,

Acad.G.Bonchev st., bl.8,

Sofia, Bulgaria

\end{document}